\documentclass[10pt]{article}
\usepackage{latexsym,amssymb,amsmath}

\def\Q{\mathbb{Q}}
\def\Z{\mathbb{Z}}
\def\R{\mathbb{R}}
\def\C{\mathbb{C}}

\def\proof{\par\noindent{\em Proof. }}
\def\eproof{\hfill{$\Box$}\bigskip}

\def\tec{\hspace{-1.6mm}{\bf. }}
\def\ds{\dots}
\def\sus{\subset}
\def\al{\alpha}
\def\be{\beta}
\def\ga{\gamma}
\def\de{\delta}
\def\ep{\varepsilon}
\def\ord{\mathrm{ord}}
\def\vni{\mathrm{int}}

\newtheorem{thm}{Theorem}[section]
\newtheorem{prop}[thm]{Proposition}

\newtheorem{lem}[thm]{Lemma}

\begin{document}
\title{Complex derivatives are continuous --- three self-contained proofs. Part 1}
\author{Martin Klazar}

\maketitle

\begin{abstract}
We prove in three ways the basic fact of analysis that complex derivatives are continuous. The first, classical, proof of 
Cauchy and Goursat uses integration. The second proof of Whyburn and Connell is topological and is based on 
the winding number. The third proof of Adel'son-Vel'skii and Kronrod employs graphs of real 
bivariate functions. We give short and concise presentation of the first proof (in textbooks it often spreads over dozens
of pages). In the second and third proof we on the contrary fill in and expand omitted steps and auxiliary results.
This part 1 presents the first two proofs. The third proof, treated in the future part 2, takes one on a tour through 
theorems and results due to Jordan, Brouwer, Urysohn, Tietze, Vitali, Sard, Fubini, Komarevsky, Young, Cauchy, Riemann 
and others. 
\end{abstract}

\section{Introduction}

Derivative of a function $f:\;\R\to\R$ need not be continuous --- the well known example $f(x)=x^2\sin(1/x)$ for $x\ne0$
and $f(0)=0$ has derivative $f'(x)=2x\sin(1/x)-\cos(1/x)$ for $x\ne0$ and $f'(0)=0$, with a discontinuity at zero. 
Fundamental distinction between analysis in real and in complex domain is that in the latter no such
example exists.

\begin{thm}[continuity of derivatives]\tec\label{zakl_veta}
If a function $f:\;\C\to\C$ has derivative $f'(z)$ for every $z\in\C$ then $f':\;\C\to\C$ is continuous. 
\end{thm}

\noindent
As is well known, complex derivatives are even more special: if $f'$ exists then $f''$ exists, so every $f^{(n)}$
for $n=1,2,\ds$ exists, and with $a_n=f^{(n)}(0)/n!$ we even have $f(z)=\sum_{n=0}^{\infty}a_nz^n$ for every
$z\in\C$. In this survey/pedagogical/historical article, motivated by teaching of basic complex analysis,
we give three self-contained proofs of Theorem~\ref{zakl_veta}. The first proof is classical, uses complex 
integration, and we attribute it to Cauchy \cite{cauch14,cauch22} and Goursat \cite{gour}. We prove Theorem~\ref{zakl_veta}
in a self-contained manner from first principles in less than 4 pages. The second proof replaces 
integration with topological arguments based on the winding number and is due to Whyburn \cite{whyb} 
and Connell \cite{conn65}. The third proof of Adel'son-Vel'skii and Kronrod \cite{adel_kron1,adel_kron2,adel_kron} 
uses geometry of graphs of real bivariate functions. After each proof we append comments. 
We want to compare complexity and length of the proofs and thus give them in details. So black boxes like 
``let $\int_{\ga}f$ be the integral of $f$ over $\ga$'' or ``let $\mathrm{ind}(\ga,p)$ be the 
index of $\ga$ with respect to the point $p$'' have to be expanded. (See Harrison \cite{harr06,harr09} 
for formalized proofs in complex analysis.) In this part 1 of our article we treat the first two proofs.
The third proof, with details of which we presently grapple, is deferred to future part 2. 

We review some terminology and basic results; more will be introduced as needed. A {\em segment $S_{a,b}\sus\C$} is 
a straight segment joining a point $a$ to a different point $b$. A {\em rectangle $R\sus\C$} means an axis-parallel closed 
rectangle with nonzero height and width. Thus $R=\{z\in\C\;|\;\al\le\mathrm{re}(z)\le\be,\ga\le\mathrm{im}(z)\le\de\}$ 
where $\al<\be$ and $\ga<\de$ are real numbers. If $\be-\al=\de-\ga$, it is a {\em square}. We replace discs and 
circles with rectangles and their boundaries as it is not straightforward to parameterize length of a circular arc 
and we want our proofs be self-contained. $U\sus\C$ always denotes a nonempty open set. A function $f:\;U\to\C$ is 
{\em holomorphic (on $U$)} if $f'(a)\in\C$ exists for every $a\in U$. This means that for every $\ep>0$ there is a 
$\de>0$ such that $|\frac{f(z)-f(a)}{z-a}-f'(a)|<\ep$ whenever $0<|z-a|<\de$, $z\in U$. We also consider 
functions $f$ holomorphic on sets $A\sus\C$ that are not open, then $f'(a)$ exists for every $a\in A$ that is not an 
isolated point of $A$. Every holomorphic function 
is continuous. An {\em entire} function is holomorphic on $\C$. If $f(z)$ is holomorphic on $A$ then so is $f(z)/z$
on $A\backslash\{0\}$ and so on. Recall that if $f:\;A\to\C$ is continuous and $A\sus\C$ is compact then $f$ is bounded, 
in fact $|f|$ attains on $A$ its extrema, and $f$ is on $A$ uniformly continuous. Also, $f(A)$ is compact.  
Compact subsets of $\C$ are exactly the closed and bounded sets. We also work with compact sets
in other spaces besides $\C$. The normed space $\C$ is complete --- every Cauchy sequence in 
$\C$ has a limit. In particular, if $\C\supset A_1\supset A_2\supset\ds$, $A_n\ne\emptyset$, are nested closed 
sets with diameters going
to $0$ then $\bigcap_{n=1}^{\infty}A_n$ is a single point. A component of an open set $A\sus\C$ means an inclusion-wise 
maximal connected (and necessarily open) subset of $A$. An open set $A\sus\C$ is connected iff any two points of $A$ 
can be joined in $A$ by a path, or equivalently by a polygonal arc.

Sharp distinction between real and complex analysis manifests in another simple way. 
The function $f(x)=1-x^2:\;\R\to\R$ has derivative $f'(x)=-2x$ everywhere and 
satisfies $1=|f(0)|>|f(x)|$ for every $x$ with $0<|x|<1$ --- its modulus $|f|$ has at $0$ a strict local maximum. 
No such complex function exists.  

\begin{thm}[maximum modulus principle]\tec\label{max_mod_pr}
If a function $f:\;\C\to\C$ has derivative $f'(z)$ for every $z\in\C$ then $|f|$ has no strict local maximum --- for 
every $a\in\C$ and $\de>0$ there is a $z\in\C$ such that $0<|z-a|<\de$ and $|f(a)|\le|f(z)|$.
\end{thm}
In fact, in both the second and third proof we deduce Theorem~\ref{zakl_veta} from Theorem~\ref{max_mod_pr}.
Logically, Theorem~\ref{max_mod_pr} is equivalent to its particular case with $a=0$, by the change of variable 
$g(z)=f(z+a)$, and this applies to Theorem~\ref{zakl_veta} too. 

We close with another example of a discontinuous derivative, taken from a 
non-archimedean realm and adapted from the write-up \cite{eth_writeup}. Let $$K=(\Q((x)),|\cdot|)$$
be the normed field of formal Laurent series with rational coefficients. Its elements $s=s(x)\in K$ 
are formal infinite series 
$s(x)=\sum_{n\ge m}a_nx^n$ with $m\in\Z$ and $a_n\in\Q$, and the norm given by $|s|=(1/2)^{\ord(s)}$ where 
$\ord(s)$ is the least $n\in\Z$ with $a_n\ne0$ ($|0|=(1/2)^{+\infty}=0$). This norm is complete but non-archimedean 
($|s+t|\le\max(|s|,|t|)$, with equality for $|s|\ne |t|$). Now consider the function $f:\;K\to K$, defined by $f(0)=0$ 
and for nonzero $s=a_mx^m+a_{m+1}x^{m+1}+\ds$, $a_m\ne0$, by
$$
f(s)=a_mx^m+a_{m+1}x^{m+1}+\ds+a_{2|m|}x^{2|m|}\;.
$$
It is easy to see that $f$ is locally constant outside zero, so $f'(s)=0$ for $s\ne0$, but $f'(0)=1$. The 
derivative of $f$ is the impulse function that has value $1$ at $0$ and $0$ elsewhere. Among real functions this is 
not a derivative as it lacks the Darboux property (real $f'$ attains all intermediate values). 

\section{Proof by integration: Cauchy and Goursat}

To prove Theorem~\ref{zakl_veta} we first develop from scratch some complex integration.
By it we prove Cauchy's theorem on vanishing of integrals and an important result on non-vanishing.
Then we deduce Cauchy's formula (in fact two), from which Theorem~\ref{zakl_veta} and much 
more follows easily.
Cauchy's formula expresses $f(a)$ for a holomorphic function $f$ in terms of values $f(z)$ for $z$ lying far
from $a$. This is the strange and fascinating non-locality of complex analysis. 

Let $a,b\in\C$, $a\ne b$, be a pair of points and $f:\;U\to\C$ be a continuous function where $U$ contains 
the segment $S=S_{a,b}=\{\varphi(t)=a+t(b-a)\;|\;t\in[0,1]\}$ spanned by $a$ and $b$. A {\em partition $P$ of $S$} is 
a tuple of points $P=(a_0,a_1,\ds,a_k)$ on $S$, with $a_i=\varphi(t_i)$ for some numbers $0=t_0<t_1<\ds<t_k=1$. 
An {\em equipartition} has $t_i=i/k$, $i=0,1,\ds,k$. We set 
$\|P\|=\max_{1\le i\le k}|a_i-a_{i-1}|$. {\em Cauchy's sum} (corresponding to $P$ and $f$) is 
$C(P,f)=\sum_{i=1}^{k}f(a_i)(a_i-a_{i-1})$. For any sequence $P_1,P_2,\ds$ of partitions of $S$ with 
$\|P_n\|\to0$ we define the {\em integral of $f$ over $S$} as
$$
\int_{a,b}f:=\lim_{n\to\infty}C(P_n,f)\;.
$$
Below we prove that the limit exists and does not depend on the choice of $P_n$; for $P_n$ one may take equipartitions 
of $S$ with $k=n$. 

Let $R$ be a rectangle. Its vertices form a quadruple of points $a,b,c,d\in\C$, ordered counter-clockwise starting 
from the lower left corner $a$. The  (oriented) {\em boundary of $R$} is $\partial R=S_{a,b}\cup S_{b,c}\cup S_{c,d}\cup S_{d,a}$ and the {\em interior} is $\vni(R)=
R\backslash\partial R$. If $f:\;U\to\C$ is continuous and $\partial R\sus U$, we define the {\em integral of $f$ over 
$\partial R$} as
$$
\int_Rf=\int_{\partial R}f:=\int_{a,b}f+\int_{b,c}f+\int_{c,d}f+\int_{d,a}f\;.
$$

\begin{prop}[on integrals]\tec\label{integ}
Let $f,g:\;U\to\C$ be continuous functions, $S=S_{a,b}\sus U$ be a segment and $R$ be a 
rectangle with $\partial R\sus U$. The limit defining $\int_{a,b}f$ always exists and is independent of
the choice of the partitions $P_n$. $\int_{a,b}$ and $\int_R$ have the following properties.
\begin{enumerate}
\item If $\al,\be\in\C$ then $\int_{a,b}(\al f+\be g)=\al\int_{a,b}f+\be\int_{a,b}g$ and similarly for $\int_R$.
\item $|\int_{a,b}f|\le\max_{z\in S}|f(z)|\cdot|b-a|$ and $|\int_Rf|\le\max_{z\in\partial R}|f(z)|\cdot p(R)$ 
where $p(R)$ is the perimeter of $R$ ($=|b-a|+|c-b|+|d-c|+|a-d|$).
\item One has $\int_{a,b}f=-\int_{b,a}f$ and if $c\in S$, $c\ne a,b$, then $\int_{a,b}f=\int_{a,c}f+\int_{c,b}f$.
\end{enumerate}
\end{prop}
\proof
First note the bound ($P=(a_0,a_1,\ds,a_k)$ is a partition of $S$) $|C(P,f)|\le\max_{z\in S}|f(z)|\cdot|b-a|$ 
which follows from $\sum_{i=1}^k|a_i-a_{i-1}|=|a_k-a_0|=|b-a|$ --- the points $a_i$ are 
collinear. It suffices to show that for every $\ep>0$ there is a $\de>0$ such that if $P_1$ and $P_2$ are partitions of $S$
with $\|P_1\|,\|P_2\|<\de$ then $|C(P_1,f)-C(P_2,f)|<\ep$. Let first $P_1\sus P_2$, $P_1=(a_0,a_1,\ds,a_k)$ and 
$\|P_1\|<\de$. Then there are partitions $Q_i$ of $S_i=S_{a_{i-1},a_i}$, $i=1,2,\ds,k$, using points of $P_2$
such that $C(P_2,f)=\sum_{i=1}^kC(Q_i,f)$. Denoting by $h$ the function that is constantly $f(a_i)$ on $S_i$ 
we have $m_i:=|f(a_i)(a_i-a_{i-1})-C(Q_i,f)|=|C(Q_i,h-f)|\le\max_{z\in S_i}|f(a_i)-f(z)|\cdot|a_i-a_{i-1}|$. 
By the uniform continuity of $f$ on $S$, for small $\de$ we have $|f(a_i)-f(z)|<\ep/|b-a|$ for every 
$z\in S_i$ and $i=1,2,\ds,k$. Thus $m_i<\ep\frac{|a_i-a_{i-1}|}{|b-a|}$ and $|C(P_1,f)-C(P_2,f)|\le m_1+\ds+m_k<\ep$.
If $P_1$ and $P_2$ are incomparable by inclusion, we take the partition $P_3=P_1\cup P_2$ and apply the previous result 
on the pairs $P_1,P_3$ and $P_2,P_3$.

1. By the linearity of $C(P,f)$ in $f$. 

2. By the bound on $|C(P,f)|$ stated at the beginning of the proof.

3. By limit transition from the corresponding identities for Cauchy's sums, which are immediate. 
For example, if $P_1$ is a partition of $S_{a,c}$ and $P_2$ of $S_{c,b}$ then the concatenated partition $P_3=P_1P_2$ of $S_{a,b}$ 
satisfies $C(P_3,f)=C(P_1,f)+C(P_2,f)$ and $\|P_3\|=\max(\|P_1\|,\|P_2\|)$.
\eproof

\begin{prop}[Cauchy's theorem]\tec\label{cauchy_thm}
Let $f:\;U\to\C$ be a continuous function and $R\sus\C$ be a rectangle. 
\begin{enumerate}
\item If $\partial R\sus U$ and $f$ is linear then $\int_Rf=0$.
\item If $R\sus U$ (the whole $R$ lies in $U$, not just its boundary) and $f$ is holomorphic then $\int_Rf=0$.
\end{enumerate}
\end{prop}
\proof
For brevity we will omit $f$ in $\int_{\dots}f$ and write just $\int_{\dots}$.

1. Let $f(z)=\al z+\be$ and the vertices of $R$ be $a,b,c,d$. If $P$ is an equipartition of $S_{a,b}$ then 
$C(P,f)=\sum_{i=1}^k(\al a_i+\be)(a_i-a_{i-1})=\al\sum_{i=1}^k(a+i(b-a)/k)(b-a)/k+
\be(b-a)=\al(a(b-a)+(\frac{1}{2}+\frac{1}{2k})(b-a)^2)+\be(b-a)$. Hence $k\to\infty$ gives 
$\int_{a,b}=\al(b^2-a^2)/2+\be(b-a)=g(b)-g(a)$ where $g(z)=\al z^2/2+\be z$. 
So $\int_R=\int_{a,b}+\int_{b,c}+\int_{c,d}+\int_{d,a}=g(b)-g(a)+g(c)-g(b)+g(d)-g(c)+g(a)-g(d)=0$. 

2. Two segments joining the midpoints of the opposite sides of $R$ divide $R$ in four rectangles $R_1,\ds,R_4$, 
each $R_i$ being congruent to $\frac{1}{2}R$. 
We see that for some $j\in\{1,2,3,4\}$, $|\int_{R_j}|\ge|\int_R|/4$. This follows 
from the triangle inequality and the identity $\int_R=\int_{R_1}+\ds+\int_{R_4}$ which follows from the definition 
of $\int_{R_i}$ and part 3 of Proposition~\ref{integ}. Clearly, $p(R_j)=p(R)/2$. Iterating this division we obtain nested rectangles 
$R=R_0\supset R_1\supset R_2\supset\ds$ such that $|\int_{R_n}|\ge|\int_R|/4^n$ and $p(R_n)=p(R)/2^n$. We take the point 
$z_0\in U$ given by
$$
\{z_0\}=\bigcap_{n=0}^{\infty}R_n
$$
by the completeness of $\C$. Since $f'(z_0)$ exists, for any $\ep>0$ for large enough $n$ and $z\in R_n$ one has
$f(z)=f(z_0)+f'(z_0)(z-z_0)+\Delta(z)(z-z_0)$ with $|\Delta(z)|<\ep$. Let $g(z)=f(z_0)+f'(z_0)(z-z_0)$ and 
$h(z)=\Delta(z)(z-z_0)$. Then $f=g+h$ and (by part 1 of Proposition~\ref{integ}) $\int_{R_n}=\int_{R_n}g+\int_{R_n}h=
\int_{R_n}h$ because $\int_{R_n}g=0$ by part 1 of this proposition. Now $|h(z)|<\ep p(R_n)$ for $z\in R_n$. Thus (by part 2 of Proposition~\ref{integ})
$$
\frac{\left|\int_R\right|}{4^n}\le\bigg|\int_{R_n}\bigg|=\bigg|\int_{R_n} h\bigg|<\ep p(R_n)^2=\frac{\ep p(R)^2}{4^n}
$$
and $|\int_R|<\ep p(R)^2$. So $\int_Rf=0$.
\eproof

\begin{prop}[non-vanishing of $\int_R$]\tec\label{nonvan}
Let $R\sus\C$ be the square with vertices $\pm1\pm i$. Then
$\rho:=\int_R1/z\ne0$.
\end{prop}
\proof
Now $a=-1-i$, $b=1-i$, $c=1+i$, and $d=-1+i$. Let $P$ be an equipartition of $S_{a,b}$. Then
$$
C(P,1/z)=\sum_{j=1}^k\frac{(b-a)/k}{a+j(b-a)/k}=\sum_{j=1}^k\frac{2/k}{2j/k-1-i}
=\sum_{j=1}^k\frac{(2/k)(2j/k-1+i)}{(2j/k-1)^2+1}
$$
and we see that $\mathrm{im}(C(P,1/z))\ge1$. Thus $\mathrm{im}(\int_{a,b}1/z)\ge1$ and in particular $\int_{a,b}1/z\ne0$. 
Since $\frac{2/k}{a+2j/k}=\frac{2i/k}{ai+2ij/k}=
\frac{2i/k}{b+2ij/k}$ and $c-b=2i$, we have $C(P,1/z)=C(iP,1/z)$ where $iP$ is the equipartition of $S_{b,c}$ obtained 
by rotating $P$ around $0$ by $i$. Thus $\int_{b,c}1/z=\int_{a,b}1/z$, and similar arguments (extending the fraction by $-1$ 
or $-i$) show that also $\int_{c,d}1/z=\int_{a,b}1/z$ and $\int_{d,a}1/z=\int_{a,b}1/z$. Hence $\int_R1/z=
4\int_{a,b}1/z\ne0$, even $\mathrm{im}(\int_R1/z)\ge4$.
\eproof

For any compact set $X\sus\C$ we set $\C_X:=\C\backslash X$ and denote by $H_X$ the set of all holomorphic functions
$f:\;\C_X\to\C$; for $X=\{a\}$ we write just $\C_a$ and $H_a$. Let $H:=\bigcup_XH_X$, $X\sus\C$ compact. 
We define the functional
$$
\int:\;H\to\C\ \mbox{ by }\ \int f:=\int_Rf\ \mbox{ where }\ \vni(R)\supset X\ \mbox{ and }\ f\in H_X\;.
$$  

\begin{prop}\tec\label{functional} 
The mapping $\int:\;H\to\C$ is correctly defined since its values do not depend on the choice of the rectangles $R$ 
and it has the following properties.
\begin{enumerate}
\item If $\al,\be\in\C$ and $f,g\in H$ then $\int(\al f+\be g)=\al\int f+\be\int g$.
\item If $f\in H_a$, $a\in\C$, and $f$ is bounded near $a$ then $\int f=0$.
\item For $(z-a)^{-1}\in H_a$ one has $\int(z-a)^{-1}=\rho\ne0$ (see Proposition~\ref{nonvan}).
\item If $f,f_n\in H_X$, $n=1,2,\ds$, $R$ is any rectangle with $\vni(R)\supset X$, and $f_n(z)\to f(z)$ 
uniformly in $z\in\partial R$, then $\int f_n\to\int f$.
\end{enumerate}
\end{prop}
\proof
Let $f\in H_X$ and $R,S\sus\C$ be two rectangles containing $X$ in their interiors. Assume first that $S$ lies in the 
interior of $R$. Using the lines extending the sides of $S$ we divide $R$ in nine rectangles $R_1,\ds, R_9$ with $S$ being
one of them, say $S=R_5$. As in the proof of part 2 of Proposition~\ref{cauchy_thm} we have 
$\int_Rf=\int_{R_1}f+\ds+\int_{R_9}f$. For $i\ne5$ we have $\int_{R_i}f=0$ by part 2 of Proposition~\ref{cauchy_thm} 
because $R_i\sus\C_X$. So $\int_R f=\int_{R_5}f=\int_Sf$. General position of $R$ and $S$ is dealt with by 
two applications of this result, by shrinking appropriately one rectangle (we use that $R\cap S$ is again a 
rectangle). So $\int_Rf=\int_Sf$ whenever $X\sus\mathrm{int}(R)\cap\mathrm{int}(S)$. 

1. This follows from part 1 of Proposition~\ref{integ}. 

2.  By part 2 of Proposition~\ref{integ}, shrink $R$ with $a\in\vni(R)$ to $a$. 

3. This follows from Proposition~\ref{nonvan} by the shift $z\mapsto z+a$.

4. By the bound in part 2 of Proposition~\ref{integ}, $|\int_R(f-f_n)|\to0$. 
\eproof

\begin{prop}[two Cauchy's formulae]\tec
Let $f:\;\C\to\C$ be entire. Then for every $a\in\C$ we have
$$
f(a)=\frac{1}{\rho}\int\frac{f(z)}{z-a}\ \mbox{ and }\ f'(a)=\frac{1}{\rho}\int\frac{f(z)}{(z-a)^2}\;.
$$
The second formula implies the continuity of $f':\;\C\to\C$ --- Theorem~\ref{zakl_veta}.
\end{prop}
\proof
Parts 1, 3 and 2 of Proposition~\ref{functional} give the first formula:
$$
\int\frac{f(z)}{z-a}=\int\frac{f(a)}{z-a}+\int\frac{f(z)-f(a)}{z-a}=f(a)\rho+0\;.
$$
For distinct $a,b\in\C$ the first formula gives
$$
\frac{f(a)-f(b)}{a-b}=\frac{1}{\rho}\int\frac{f(z)}{(z-a)(z-b)}\;.
$$
Let $R$ be a rectangle and $a\in\vni(R)$, then $\frac{1}{(z-a)(z-b)}-\frac{1}{(z-a)^2}=O(a-b)$ if $b\in\vni(R)$, 
$z\in\partial R$ and $b$ stays away from $\partial R$. Thus if $b_n\to a$, $b_n\in\mathrm{int}(R)$ and $b_n\ne a$, 
then $\frac{f(z)}{(z-a)(z-b_n)}\to\frac{f(z)}{(z-a)^2}$ uniformly in $z\in\partial R$ (here $X=\{a,b_1,b_2,\ds\}$ and   
$f(z)$ is bounded on $\partial R$). Part 4 of Proposition~\ref{functional} gives the second formula for $f'(a)$. 
The continuity of $f'$ follows in much the same way because if $R$ is any rectangle then $\frac{1}{(z-a)^2}-
\frac{1}{(z-b)^2}=O(a-b)$ if 
$a,b\in\vni(R)$, $z\in\partial R$ and both $a,b$ stay away from $\partial R$. So if $b_n\to a$ then 
$f'(b_n)\to f'(a)$ by the second formula and part 4 of Proposition~\ref{functional}.
\eproof

\noindent
This concludes the first proof of Theorem~\ref{zakl_veta}.

\subsection{Comments}

Next we easily obtain a power series expansion $f(a)=\sum_{n\ge0}a_na^n$: use part 4 of 
Proposition~\ref{functional} with $X=\{0,a\}$ and $f_n(z)=f(z)\sum_{k=0}^n\frac{a^k}{z^{k+1}}$. By the identity
$1+y+y^2+\ds+y^n=\frac{1-y^{n+1}}{1-y}$, $f_n(z)\to\frac{f(z)}{z-a}$ uniformly in $z\in\partial R$ for any 
rectangle $R$ with $0,a\in\vni(R)$ and the distance of $\partial R$ and $a$ exceeding $|a|$. The power series 
expansion with coefficients $a_n=\frac{1}{\rho}\int f(z)/z^{n+1}$ now follows by Cauchy's first formula. With a
little more work one shows that each $f^{(n)}$ exists and $a_n=f^{(n)}(0)/n!$.

The above proof of Theorem~\ref{zakl_veta}, based on complex integration and Cauchy's formula, is classical and 
appears in some form in every textbook on complex analysis. Initially we just wanted to 
give a short self-contained presentation of it. For example, Ahlfors \cite{ahlf} 
reaches Cauchy's formula by page 95, Bak and Newman \cite{bak_newm} by page 58, \v Cern\'y \cite{cern} by page 202, 
Conway \cite{conw} by page 84, Henrici \cite{henr} by page 245, Krantz \cite{kran} by page 26 (but this is 
a handbook without proofs), Lin \cite{lin} by page 394, Markushevich \cite{mark} by page 162, Palka  by page 161, 
Rudin \cite{rudi} by page 207, Titchmarsh \cite{titc} by page 99, Vesel\'y \cite{vese} by page 105, and so on, 
but certainly one should not need 26 or 207 or so many pages to prove it? To be fair, these works do not set as their goal
minimalistic proof of Cauchy's formula or Theorem~\ref{zakl_veta}, but the question how such proof would look like is 
legitimate and natural. The question if there is an integration-free proof of ``$f'$ exists $\Rightarrow f''$ exists'' 
for complex functions was posed at mathoverflow \cite{maov} in 2011, and the discussion there pointed to the Whyburn and Connell
type proof which we present in the next section. Literature did not give us an answer and we 
wrote such minimalistic proof ourselves above. We attribute it to Cauchy and Goursat because of their techniques 
and ideas behind it but the actual concise presentation is of course ours. 

Goursat \cite{gour} was the first to prove Cauchy's theorem (Proposition~\ref{cauchy_thm} and its variants) assuming 
only the existence of $f'(z)$ and not its continuity (Cauchy assumed continuity of the derivative 
in his arguments). See Zalcman \cite{zalc} for further information on proofs of Cauchy's theorem. Cauchy's sums 
in the definition of $\int_{a,b}f$ are indeed due to Cauchy and he used them to introduce real 
and complex integrals, see \cite[lectures 21 and 23]{cauch23} (we learned it in Schwabik and \v Sarmanov\'a 
\cite{schw_sarm}). For the information on when, where and how actually Cauchy proved 
his eponymous theorem and formula we turned to the thorough historical work of Smithies \cite{smit}. Smithies states that Cauchy's 
theorem appears implicitly already in Clairaut \cite{clai} in 1743 and that also Euler mentions it implicitly 
(\cite[p. 7]{smit}) and that we find it first in Cauchy's work, for rectangles and a class of curvilinear quadrilaterals,  
in the 1814 memoir \cite{cauch14} that was published 13 years later (\cite[p. 56]{smit}). Cauchy's formula appears
in a certain form in Cauchy's 1822 article \cite{cauch22} but ``it was not until about 1831 (cf. Section 6.9 below) 
that Cauchy perceived how powerful a tool this equation could be.'' (\cite[p. 71 and 84]{smit}).

Proposition~\ref{nonvan} and its proof is due to the author. We could compute the well-known value 
$\rho=\int 1/z=2\pi i$ but this is not necessary to prove Theorem~\ref{zakl_veta}. Cauchy's theorem is often and 
rightly hailed as a cornerstone of complex analysis, but the non-vanishing in Proposition~\ref{nonvan} is equally 
important: if everything vanished there would be no Cauchy's formula.  

\section{Proof by the winding number: Whyburn and Connell}

We prove Theorem~\ref{zakl_veta} in four steps. We introduce the winding number. We deduce by it the maximum modulus principle, 
which is the hardest step. From the maximum modulus principle we deduce Theorem~\ref{zakl_veta}. We return 
to the winding number and prove its existence. The proof avoids integration.

Let $I=[0,1]$. A {\em path} is a continuous mapping $f:\;I\to\C$. An {\em arc} is an injective path. 
A {\em loop} is a path such that $f(0)=f(1)$.
For a loop $f$ and a point $p\in\C\backslash f(I)$ the {\em winding number $w(f,p)\in\Z$} counts how many times $f(t)$ 
winds around $p$ when $t$ moves in $I$ from $0$ to $1$. Proposition~\ref{wind_numb} below summarizes 
properties we need $w(f,p)$ to have; its proof is deferred to the end. A {\em circuit $f$} is  a 
loop that is injective with the single exception $f(0)=f(1)$. In a {\em positive circuit $f$}  
for $t\in I$ moving from $0$ to $1$ the interior of $f(I)$ (the bounded component of $\C\backslash f(I)$) lies to the 
left of $f(t)$. This notion is not easy to make rigorous in general as it is already hard to prove 
existence of the interior for a general circuit, see Hales \cite{hale1, hale2} and part 2 of this article. 
But we will use it only for rectangular circuits and then positivity and interiors are clear. (Unlike the third proof, 
we do not need general Jordan curve theorem.) We say that two loops $f,g$ are {\em homotopic 
in $X\sus\C$} if there is a continuous mapping $h:\;I\times I\to X$, called a {\em homotopy}, such that 
for every $t,u\in I$ one has $h(0,u)=h(1,u)$, $h(t,0)=f(t)$ and $h(t,1)=g(t)$. We also say that $f$ {\em can be 
deformed to $g$ in $X$}. Two loops $f,g:\;I\to\C$ are 
{\em equivalent} if they differ only in starting points, which means that for an $s\in I$ one has 
$g(t)=f(s+t)$ for every $t\in I$ with $s+t$ taken modulo $1$. A {\em product $h=f*g$} of two paths $f$ and $g$ with $f(1)=g(0)$ 
is the path $h$ given by $h(t)=f(2t)$ for $t\in[0,\frac{1}{2}]$ and $h(t)=g(2t-1)$ for $t\in[\frac{1}{2},1]$.

\begin{prop}[winding number]\tec\label{wind_numb}
There exists a mapping
$$
w:\;\{(f,p)\;|\;\mbox{$f$ is a loop},\; p\in\C\backslash f(I)\}\to\Z
$$
with the following properties.
\begin{enumerate}
\item If $f$ is constant then $w(f,p)=0$. If $f$ is a positive circuit such that 
$f(I)=\partial R$ for a rectangle $R$ and $p\in\vni(R)$ then $w(f,p)=1$. 
\item If $f$ is a loop and  $p,q\in\C\backslash f(I)$ are points lying in the same component of $\C\backslash f(I)$ 
then $w(f,p)=w(f,q)$. 
\item If $p$ is a point and $f,g$ are loops that are homotopic in $\C\backslash\{p\}$ then $w(f,p)=w(g,p)$. 
\item If $f$ and $g$ are two equivalent loops and $p\in\C\backslash f(I)=\C\backslash g(I)$ then $w(f,p)=w(g,p)$.
\item If $h=f*g$ is a product of loops and $p\in\C\backslash h(I)$ then $w(h,p)=w(f,p)+w(g,p)$.
\end{enumerate}
\end{prop}
We prove it at the end of the section. Now we use $w(f,p)$ to get the maximum modulus principle. 
One obtains it easily by Cauchy's formula or by local power series expansion. But if one wants to avoid integration 
and cannot use these tools, things are not so easy.

For a rectangle $R$ we denote by $\varphi_R:\;I\to\partial R$ the circuit traversing $\partial R$
once counter-clockwisely. More precisely, if the vertices of $R$ are $(a,b,c,d)=(a_0,\ds,a_3)$
we define $\varphi_R(t)=a_i+(a_{i+1}-a_i)(4t-i)$ for $\frac{i}{4}\le t\le\frac{i+1}{4}$ where $i=0,1,2,3$ and $a_4=a_0=a$.
If $f:\;\partial R\to\C$ is continuous, we denote by $f_R$ the loop $f_R=f\circ\varphi_R$. 

\begin{prop}[properties of $w(f,p)$]\tec\label{cor_w_num}
The winding number has the following properties.
\begin{enumerate}
\item If $R$ is a rectangle, $f:\;R\to\C$ is continuous, $p\in\C\backslash f(\partial R)$, 
$W$ is the component of $\C\backslash f(\partial R)$ containing $p$, and $w(f_R,p)\ne0$, then 
$$
f(\vni(R))\supset W\;.
$$ 
\item If $R$ is a rectangle, $f:\;R\to\C$ is continuous, and $z_0\in\vni(R)$ is a point such that $f'(z_0)\ne0$, then 
for any sufficiently small square $S\sus R$ centered at $z_0$ one has 
$$
w(f_S,f(z_0))=1\;.
$$
\item If $R$ is a rectangle, $f:\;R\to\C$ is continuous, and $p\in\C\backslash f(\partial R)$ is 
a point such that $f^{-1}(p)=\{z_1,z_2,\ds,z_k\}\sus\vni(R)$ (possibly $k=0$), then for any sufficiently small 
disjoint squares $S_i\sus\vni(R)$, $i=1,2,\ds,k$ and $S_i$ is centered at $z_i$, one has 
$$
w(f_R,p)=\sum_{i=1}^kw(f_{S_i},p)\;.
$$ 
\end{enumerate}
\end{prop}
\proof
1. Let $q\in W$ be an arbitrary point. Then $w(f_R,q)=w(f_R,p)\ne0$ by part 2 of Proposition~\ref{wind_numb}.
If $q\not\in f(\vni(R))$ then by parts 1 and 3 of Proposition~\ref{wind_numb} we get a contradictory value 
$w(f_R,q)=0$ as then $f_R$ is homotopic in $\C\backslash\{q\}$ to a constant mapping (we continuously shrink $R$ to an 
interior point and take restrictions of $f$ to the shrunk rectangles). So $q\in f(\vni(R))$. 

2. Let $g(z)=f(z_0)+f'(z_0)(z-z_0)$. For $S$ small enough and $z=\varphi_S(t)$, 
$$
|f_S(t)-g_S(t)|=o(z-z_0)<|g_S(t)-f(z_0)|=|f'(z_0)(z-z_0)|
$$ 
for every $t\in I$. By parts 3 and 1 of Proposition~\ref{wind_numb}, 
$$
w(f_S,f(z_0))=w(g_S,f(z_0))=1\;. 
$$
The first equality sign is due to the homotopy of $f_S$ and $g_S$ in 
$\C\backslash\{f(z_0)\}$ via $h(t,u)=f_S(t)+u(g_S(t)-f_S(t))$. By the displayed inequality, $h(t,u)$
avoids the point $f(z_0)$ because the segment spanned by $f_S(t)$ and $g_S(t)$ is too short to contain it. 
The second equality sign follows from the second claim in part 1 of Proposition~\ref{wind_numb} and by 
the homotopy of $g_S$ and $j_S$, $j(z)=f(z_0)+z-z_0$, in $\C\backslash\{f(z_0)\}$ via $h(t,u)=
f(z_0)+s(u)(\varphi_S(t)-z_0)$, where $s(u)$ is any arc joining the points $f'(z_0)$ and $1$ and avoiding $0$ 
(then $h(t,u)$ clearly avoids $f(z_0)$).

3. Perturbing $f$ a little (and using part 3 of Proposition~\ref{wind_numb}) we may assume that 
$\mathrm{re}(z_1)<\mathrm{re}(z_2)<\ds<\mathrm{re}(z_k)$. We take small enough squares $S_i\sus\vni(R)$
centered at $z_i$ such that $\mathrm{re}(S_1)<\mathrm{re}(S_2)<\ds<\mathrm{re}(S_k)$. Let $a,b,c,d$ be the 
vertices of $R$ and $a_i,b_i,c_i,d_i$ be those of $S_i$, $i=1,2,\ds,k$. By $L_i$ we denote the vertical segment
joining $a_i$ to a point $e_i\in S_{a,b}$. We define a loop $\psi:\;I\to T=\partial R\cup\bigcup_{i=1}^k(\partial S_i
\cup L_i)$ as follows. We start at $a$ and traverse linearly
in this order the segments $S_{a,e_1},S_{e_1,a_1}$, $\partial S_1$ clockwise, $L_1=S_{a_1,e_1}$, $S_{e_1,e_2}$, 
$S_{e_2,a_2}$, $\partial S_2$ clockwise, $L_2=S_{a_2,e_2}$, $\ds$, $L_k=S_{a_k,e_k}$, $S_{e_k,b}$, $S_{b,c}$, $S_{c,d}$,
and $S_{d,a}$. Thus each $L_i$ is traversed twice, first up and then down, but otherwise $\psi$ is injective. We denote 
$f_T=f\circ\psi$ and show that 
$$
w(f_R,p)-\sum_{i=1}^kw(f_{S_i,p})=w(f_T,p)=0\;,
$$
which will prove the claim. 

The first equality follows by repeated applications of parts 3, 4, 5, and 1 of 
Proposition~\ref{wind_numb}. Using part 4 we move the starting point of $\psi$ from $a$ to $e_1$ (more precisely, we
appropriately deform the time schedule of $\psi$ so that it visits $e_1$ at the times $t=0,\frac{1}{2},1$). 
We denote by $T_1$ the initial loop of $\psi$ going from $e_1$ to $e_1$ via $L_1$ and $\partial S_1$, 
and by $U_1$ the rest of $\psi$. Part 5 then gives $w(f_T,p)=w(f_{T_1},p)+w(f_{U_1},p)$. Applying parts 4 and 5 to $T_1$, 
we see that $w(f_{T_1},p)=-w(f_{S_1},p)$ because traversing $L_1$ back and forth gives (by parts 3 and 1) winding $0$,
and $\partial S_1$ is traversed clockwise. (Using parts 5, 3, and 1 one easily shows that reverting direction of 
a loop changes the sign of the winding number because the product of a loop and its reversal is homotopic outside 
the point to a constant.)
So $w(f_T,p)=-w(f_{S_1},p)+w(f_{U_1},p)$. Then we move the starting point of $U_1$ to $e_2$ and argue similarly
(note that $U_1$ is $\psi$ without the loop $T_1$). After $k$ steps we get $w(f_T,p)=-w(f_{S_1},p)-\ds-w(f_{S_k},p)+w(f_{U_k},p)$. 
But $U_k$ traverses counter-clockwisely $\partial R$, only starts and finishes at $e_k$ 
instead of $a$, so $w(f_{U_k},p)=w(f_R,p)$ by part 4.

The second equality follows from parts 3 and 1 of Proposition~\ref{wind_numb} by deforming $f_T$ to a constant mapping 
in $\C\backslash\{p\}$. We move each segment $L_i$ continuously 
to the right to the vertical segment $S_{b_i,f_i}$, $f_i\in S_{a,b}$, so that the parts $S_{b_i,a_i},S_{a_i,e_i},S_{e_i,f_i}$ of 
$\psi$ become $S_{b_i,f_i}$. Then we move the segments of $\psi$ lying on $S_{a,b}$, which are $S_{a,e_1},S_{f_1,e_2},S_{f_2,e_3},\ds,
S_{f_k,b}$, continuously upward and take with them the segments $S_{d_i,c_i}$ as encountered, to the position with 
the imaginary part $v=\max_{1\le i\le k}\mathrm{im}(d_i)$. We take compositions $f\circ\psi$ with the closed 
broken line $\psi$ being deformed and accomplish in this way the deformation of $f_T$ in $\C\backslash\{p\}$ to 
$f_U$ where $U$ is the subrectangle $U=\{z\in R\;|\;\mathrm{im}(z)\ge v\}$. By shrinking $U$ continuously 
to an interior point we complete the deformation of $f_T$ in $\C\backslash\{p\}$ to a constant.
\eproof

For a rectangle $R$ with vertices $a,b,c,d$ we denote by $|R|$ its {\em area}, $|R|=|a-b|\cdot|c-b|>0$.
If $R$ is given by $\mathrm{re}(z)\in[\al,\be]$ and $\mathrm{im}(z)\in[\ga,\de]$ then 
a {\em grid partition of $R$} is given by some partitions $\al=t_0<t_1<\ds<t_k=\be$ and $\ga=u_0<u_1<\ds<u_l=\ga$ 
and consists of the rectangles $S_1,S_2,\ds,S_{kl}$ given by the intervals $[t_{i-1},t_i]$ and $[u_{j-1},u_j]$. 
Clearly, the $S_n$ have disjoint interiors, $R=\bigcup_{n=1}^{kl}S_n$, and $|R|=\sum_{n=1}^{kl}|S_n|$. 
Recall that a set $A\sus\C$ {\em has measure $0$} if for every $\ep>0$ there is a finite or infinite sequence of 
rectangles $R_1,R_2,\ds$ such that $A\sus R_1\cup R_2\cup\ds$ and $|R_1|+|R_2|+\ds<\ep$. Clearly, every subset 
of measure $0$ set has measure $0$ and every countable union of measure $0$ sets has measure $0$.
For $A\sus\C$ we denote by $d(A)=\sup_{a,b\in A}|a-b|$ the {\em diameter of $A$.}

\begin{prop}[on area]\tec\label{area}
Let $R$ be a rectangle. The area function $|R|$ has the following properties.
\begin{enumerate}
\item If $R_n\sus R$, $n=1,2,\ds$, is a finite or infinite sequence of rectangles with disjoint interiors then 
$|R_1|+|R_2|+\ds\le|R|$.
\item If $R_n$, $n=1,2,\ds$, is a finite or infinite sequence of rectangles with
$R_1\cup R_2\cup\ds\supset R$ then $|R_1|+|R_2|+\ds\ge|R|$. Hence no rectangle has measure $0$ and 
$R\backslash A\ne\emptyset$ for every measure $0$ set $A$.
\end{enumerate}
\end{prop}
\proof
1. It clearly suffices to consider only a finite sequence $R_1,R_2,\ds,R_n$ of such rectangles. We divide $R$ by the 
lines extending the sides of $R_i$ into a grid partition $S_1,S_2,\ds,S_m$, so 
$|S_1|+\ds+|S_m|=|R|$. For $i=1,2,\ds,n$ there are subsets $X_i\sus\{1,2,\ds,m\}$, pairwise disjoint, such that 
$R_i=\bigcup_{j\in X_i}S_j$. Again, $|R_i|=\sum_{j\in X_i}|S_j|$ since the $S_j$ form a grid partition of $R_i$. 
Thus $|R_1|+\ds+|R_n|$ equals to a subsum of $|S_1|+\ds+|S_m|=|R|$ and the inequality follows.

2. Again we need to consider only a finite sequence $R_1,R_2,\ds,R_n$ of such rectangles. The reason
is compactness of $R$: If we magnify $R_1,R_2,\ds$ by a factor $c>1$ to $R_1',R_2',\ds$ then the sum of areas 
magnifies by $c^2$ and $\vni(R_1')$, $\vni(R_2'),\ds$ is an open covering of $R$. It has a finite subcovering given by the
set of indices $J$. If we show that $\sum_{n\in J}|R_n'|\ge|R|$ we are done since $\sum_n|R_n|\ge 
c^{-2}\sum_{n\in J}|R_n'|$ and $c$ may be as close to $1$ as we wish. (Without magnifying, $R_1,R_2,\ds$ need not 
possess finite subcovering of $R$.) Using the lines extending the sides of 
$R_1,\ds,R_n$ and of $R$ we get again rectangles $S_1,S_2,\ds,S_m$ that have disjoint interiors and such that for 
$i=0,1,\ds,n$ there are subsets $X_i\sus\{1,2,\ds,m\}$, this time typically intersecting, with 
$R_i=\bigcup_{j\in X_i}S_j$, $R_0=R$, and $S_j$ forming a grid partition of $R_i$. 
Again $|R_i|=\sum_{j\in X_i}|S_j|$ for $i=0,1,\ds,n$. Moreover, the covering assumption means that 
$X_0\sus X_1\cup\ds\cup X_n$. Thus $|R|=\sum_{j\in X_0}|S_j|$ equals to a subsum of $|S_1|+\ds+|S_m|=|R_1|+\ds+|R_n|$ 
and the inequality follows.
\eproof

A {\em rooted tree} is a triple $T=(r,T,p)$ where $T$ is a set of {\em vertices}, $r\in T$ is a {\em root}, and 
$p:\;T\backslash\{r\}\to T$, the {\em child to parent mapping}, has the property that for every $v\in T\backslash\{r\}$, 
after finitely many applications of $p$ one has $p(p(\ds(p(v))\ds))=r$. {$T$ \em finitely branches} if 
every $p^{-1}(v)$, $v\in T$, is finite. A {\em branch} in $T$ is any sequence of vertices $v_1,v_2,\ds$ such that 
$v_1=r$ and $p(v_{n+1})=v_n$ for every $n$. 

\begin{lem}[K\"onig's]\tec\label{truekonig}
Every infinite finitely branching rooted tree $T=(r,T,p)$ has an infinite branch.
\end{lem}
\proof
We call $v\in T$ a {\em descendant of $u\in T$} if $p(p(\ds(p(v))\ds))=u$ after some applications of $p$. We define $v_1=r$
and $v_2\in p^{-1}(v_1)$ to be the child of $v_1$ with infinitely many descendants. Such a child exists because any 
finite union of finite sets is finite. We set $v_3$ to be the child of $v_2$ with infinitely many descendants and 
so on. This gives an infinite branch $v_1,v_2,v_3,\ds\;$.
\eproof

\begin{prop}[covering lemma]\tec\label{konig}
Let $R$ be a square, $A\sus R$ be a nonempty subset, and $F$ be a family of squares in $R$ satisfying that 
for every $a\in A$ there is a $\de>0$ such that any square $S\sus R$ with $a\in S$ and $d(S)<\de$ belongs to $F$.
Then the following holds.
\begin{enumerate}
\item A subfamily $G\sus F$ exists that is at most countable, its members have disjoint interiors, and 
$\bigcup G\supset A$.
\item If $A$ is closed then the subfamily $G$ in part 1 can be chosen to be finite. 
\end{enumerate}
Analogous results hold for segments in place of squares. 
\end{prop}
\proof 
Let $T=(r,T,p)$ be the infinite rooted $4$-ary tree where $r=R$, 
$$
T=\{R\}\cup\{\mbox{the quarters of $R$}\}\cup\{\mbox{the quarters of the quarters of $R$}\}\cup\ds
$$
(quarters of a square arise by cutting it with the two segments joining the midpoints of the opposite sides)
and $p$ sends each quarter to its parent square. 
Clearly, $\bigcup p^{-1}(S)=S$ for every $S\in T$. We view $T$ also as a poset $(T,\sus)$. Note that $S_1\in T$ 
is a descendant of $S_2\in T$ 
iff $S_1\sus S_2$ and that if $S_1$ and $S_2$ are incomparable by $\sus$ then they have disjoint interiors. 
Every strictly ascending chain in $(T,\sus)$ is finite, and so for any set of vertices $U\sus T$ every $S\in U$ is 
contained in a unique maximal element of $U$.  

1. Take $U=\{S\in T\cap F\;|\;S\cap A\ne\emptyset\}$ and set $G$ to be the maximal elements 
of $U$. It is easy to see that $G$ has the stated properties. 

2. We assume that $R\not\in F$ (else we set $G=\{R\}$) and define $U\sus T$ as consisting of 
$R$, of all children of $R$ that intersect $A$ and are 
not in $F$, of all children of these children that intersect $A$ and are not in $F$, and so on. We claim that 
the subtree $U=(r,U,p\,|\,U)$ of $T$ is finite. If not, the infinite branch $R=S_1\supset S_2\supset\ds$
in $U$ provided by Lemma~\ref{truekonig} would produce a point $\{z_0\}=\bigcap_{n=1}^{\infty}S_n$. But $z_0\in A$ 
because $A$ is closed, and so for large $n$ (since $d(S_n)\to0$) the assumption on $F$ implies that $S_n\in F$, 
which contradicts the definition of $U$. So $U$ is finite. We define $G$ to be all those children of the vertices in 
$U$ that lie in $T\backslash U$ and intersect $A$. It is not hard to see that $G\sus F$, $G$ is finite, its members 
are pairwise incomparable by $\sus$, and $\bigcup G\supset A$.

The version for segments is clear: take the infinite rooted binary tree of segments obtained by repeated cutting 
into two halves (in midpoints).
\eproof

Let $\emptyset\ne X\sus\C$ be any set, possibly not open. We call more generally a function $f:\;X\to\C$ {\em holomorphic} 
if it has derivative at every point of $X$ that is not isolated (derivative is not well defined at an isolated point). 

\begin{prop}[zero derivatives]\tec\label{zero_deri}
We have the following results on vanishing of derivatives. 
\begin{enumerate}
\item Let $X\sus\C$ be a nonempty compact set, $f:\;X\to\C$ be holomorphic, and $p\in\C$ be a point such that 
$f^{-1}(p)=\{z\in X\;|\;f(z)=p\}$ is infinite. Then $f'(b)=0$ for some $b\in f^{-1}(p)$. 
\item Let $L\sus\C$ be a segment and $f:\;L\to\C$ be a function such that $f'(z)=0$ for every $z\in L$. Then $f(L)$ is
a single point.
\item Let $A\sus\C$ be bounded and $f:\;A\to\C$ be a function such that $f'(a)=0$ for every non-isolated point 
$a\in A$. Then $f(A)$ has measure $0$.
\end{enumerate}
\end{prop}
\proof
1. The set $f^{-1}(p)$ is bounded and closed in $X$, hence is closed in $\C$ and therefore is compact. Since it is infinite, 
it has a limit point: there exist points $b_n,b\in f^{-1}(p)$ such that $b_n\ne b$ but $b_n\to b$, so $b$ is not 
isolated. But then $\frac{f(b)-f(b_n)}{b-b_n}=\frac{p-p}{b-b_n}=0$ for every $n$ and $f'(b)=0$. 

2. We fix a small number $c>0$ and consider the family $F$ of subsegments $K\sus L$ for which there 
is a $z_0\in K$ such that for every $z\in K$ one has $f(z)-f(z_0)=f'(z_0)(z-z_0)+\Delta(z)(z-z_0)=\Delta(z)(z-z_0)$ 
with $|\Delta(z)|\le c$. By part 2 of 
Proposition~\ref{konig} applied to the segment $L$ and $A=L$, there is a partition $P=(a_0,a_1,\ds,a_k)$ of $L$
such that each segment $S_{a_{i-1},a_i}\in F$. It follows that $|f(a_{i-1})-f(a_i)|\le c|a_{i-1}-a_i|$ and 
$|f(a_k)-f(a_0)|=|\sum_{i=1}^k(f(a_i)-f(a_{i-1}))|\le c\sum_{i=1}^k|a_i-a_{i-1}|=c|a_k-a_0|$. We are done 
--- $c$ may be arbitrarily small, hence $f(a_k)=f(a_0)$, and $a_0,a_k$ may be any two points of $L$, not just 
its endpoints. 

3. We may assume that $A$ has no isolated point and $f'=0$ everywhere on $A$ (since the set of isolated points is 
countable). We fix a small number $c>0$ and a square $R\supset A$. Let $F$ be the family of squares $S\sus R$ 
such that for a $z_0\in S\cap A$ and every $z\in S\cap A$,
$$
f(z)-f(z_0)=f'(z_0)(z-z_0)+\Delta(z)(z-z_0)=\Delta(z)(z-z_0)\ \mbox{ where }\ |\Delta(z)|<c\;.
$$ 
Hence $d(f(S\cap A))\le2c\cdot d(S\cap A)$. We apply to $R$, $A$, and the family $F$ part 1 of 
Proposition~\ref{konig} and cover $A$ by squares $R_1,R_2,\ds$ in $F$ with disjoint interiors. Every bounded set 
$B\sus\C$ is contained in a square $S$ with $d(S)\le 2\sqrt{2}d(B)\le 3d(B)$. For every square $S$ we have 
$|S|=d(S)^2/2$. Thus there are squares $S_n$ such that $f(R_n\cap A)\sus S_n$ --- they cover $f(A)$ --- and 
$|S_n|=d(S_n)^2/2\le(9/2)d(f(R_n\cap A))^2\le18c^2d(R_n)^2=36c^2|R_n|$. Since the $R_n\sus R$ have disjoint interiors, 
by part 1 of Proposition~\ref{area} we have $\sum_n|S_n|\le36c^2\sum_n|R_n|\le36c^2|R|$: $f(A)$ has measure $0$.
\eproof

\begin{prop}[maximum modulus principle]\tec\label{mmodp}
Let $R$ be a rectangle and $f:\;R\to\C$ be a continuous function that is holomorphic on $R\backslash P$ where $P$ 
is a finite set. Then 
$$
\max_{z\in R}|f(z)|=\max_{z\in\partial R}|f(z)|
$$
--- $|f|$  attains its maximum always on the boundary.
\end{prop}
\proof
Let $z_0\in\vni(R)$. If $f(z_0)\in f(\partial R)$ then $|f(z_0)|\le\max_{z\in\partial R}|f(z)|$.
We assume that $f(z_0)\not\in f(\partial R)$ and show that $\max_{z\in R}|f(z)|$ is not attained at $z=z_0$.
Let $W$ be the component of $\C\backslash f(\partial R)$ containing $f(z_0)$. We take a point $z_1\in\vni(R)$
near $z_0$ and a segment $L\sus R\backslash P$ such that $f(z_1)\in W$ and $L$ joins $z_1$ with a point in $\partial R$ 
(if $z_0\not\in P$ we may take $z_1=z_0$). Then $f(L)$ joins $f(z_1)$ with a point in $f(\partial R)$. 
Since $f(L)\ne\{f(z_1)\}$ and $W$ is open, there is a point $u\in L$ such that $f(S_{z_1,u})\sus W$ and $f(u)\ne f(z_1)$. 
By part 2 of Proposition~\ref{zero_deri}, $f'\ne0$ somewhere on $S_{z_1,u}$: we get a point $z_2\in\vni(R)$ such 
that $f'(z_2)\ne0$ and $f(z_2)\in W$. By part 2 of Proposition~\ref{cor_w_num}, there is a square 
$S\sus R$ centered at $z_2$ such that $f(S)\sus W$ and $w(f_S,f(z_2))=1\ne0$. Let $V$ be the component of 
$\C\backslash f(\partial S)$ containing $f(z_2)$. By part 1 of Proposition~\ref{cor_w_num}, $f(\vni(S))\supset V$. 
Let $T$, $T\sus V\sus W$, be a rectangle and $A=\{z\in R\;|\;f'(z)=0\}$. By part 3 of 
Proposition~\ref{zero_deri} and part 2 of Proposition~\ref{area}, $T\backslash f(A\cup P)\ne\emptyset$ because 
$f(A\cup P)$ has measure $0$. We select a point $z_3\in T\backslash f(A\cup P)$. Then $f^{-1}(z_3)\ne\emptyset$ as 
$z_3\in f(\vni(S))$, $f^{-1}(z_3)\sus\vni(R)$ as $z_3\in W$ and $W$ is disjoint to $f(\partial R)$, and 
$f'$ exists and is nonzero on $f^{-1}(z_3)$ as $z_3\not\in f(A\cup P)$. By part 1 of Proposition~\ref{zero_deri}, 
$f^{-1}(z_3)$ is a finite set with $k\ge1$ elements. By part 2 of Proposition~\ref{cor_w_num}, for any small enough 
and mutually disjoint squares $S_i\sus\vni(R)$, $i=1,2,\ds,k$, centered at the elements of $f^{-1}(z_3)$ one has
$w(f_{S_i},z_3)=1$. By part 3 of Proposition~\ref{cor_w_num}, $w(f_R,z_3)=\sum_{i=1}^kw(f_{S_i},z_3)=k\ne0$. By part 1 
of Proposition~\ref{cor_w_num}, $f(\vni(R))\supset W$. We are done: $f(z_0)\in W$, $W$ is open, and $f(\vni(R))\supset W$, 
thus $f$ has values with modulus larger than $|f(z_0)|$. 
\eproof

To deduce continuity of derivatives is now easy. In the following we define for $z=x$ the fraction 
$\frac{f(z)-f(x)}{z-x}$ as $f'(x)$, if this derivative exists. 

\begin{prop}[derivative is continuous]\tec\label{end_2nd_proof}
Let $R\sus\C$ be a rectangle and $f:\;R\to\C$ be a holomorphic function. 
Then for every $\ep>0$ and every rectangle $S\sus\vni(R)$ there is a $\de>0$ such that for every 
$z,x\in S$ with $|z-x|<\de$ we have
$$
\bigg|f'(z)-\frac{f(z)-f(x)}{z-x}\bigg|<\ep\;.
$$
Hence $f':\;\vni(R)\to\C$ is continuous, which gives Theorem~\ref{zakl_veta}. 
\end{prop}
\proof
The set $\partial R\times S$ is compact, hence the function $(y,x)\mapsto\frac{f(y)-f(x)}{y-x}$ on it is 
uniformly continuous, and for given $\ep>0$ there is a $\de>0$ 
such that for every $y\in\partial R$ and $x,z\in S$ with $|x-z|<\de$ we have
$$
\bigg|\frac{f(y)-f(x)}{y-x}-\frac{f(y)-f(z)}{y-z}\bigg|<\ep\;.
$$
For fixed $x,z\in S$ with  $|x-z|<\de$ we denote by $h(y):\;R\to\C$ the function given by the above difference
in the absolute value. It is continuous on $R$, holomorphic on $R\backslash\{x,z\}$, and $|h|<\ep$ on $\partial R$. 
By Proposition~\ref{mmodp}, $|h|<\ep$ on $R$. For $y=z$ we obtain the stated inequality. 

To show continuity of $f'$, fix $z,x\in\vni(S)$ with $x\ne z$ and $|z-x|<\de$. For any $u\in\vni(S)$ 
with sufficiently small $|z-u|$ is $|\frac{f(z)-f(x)}{z-x}-\frac{f(u)-f(x)}{u-x}|$ as small as we need. 
The stated inequality and the triangle inequality then imply $|f'(z)-f'(u)|<2\ep$ and we see that $f'$ is continuous
at $z$.
\eproof

But it remains to prove Proposition~\ref{wind_numb}. Let a loop $f:\;I\to\C$ and a point $p\in\C\backslash f(I)$ 
be given. We take the square $R=R_p$ centered at $p$ and with perimeter $1$ (the lower left corner of $R$ is 
$a=p-\frac{1}{8}-\frac{i}{8}$). Let $r,s\in\partial R$ be a pair of points. By $l(r,s)\in[0,\frac{1}{2}]$ 
we denote the length of the shorter of the two arcs of $\partial R$ joining $r$ and $s$ (so $l(r,s)=0$ 
if $r=s$ and maximally 
$l(r,s)=\frac{1}{2}$ when the two arcs have equal lengths), and we set $\ell(r,s)=\pm l(r,s)$ with the sign 
$+$ if $l(r,s)=\frac{1}{2}$ or if $l(r,s)<\frac{1}{2}$ and the shorter arc goes counter-clockwise in the 
$r\to s$ direction, and $-$ if it goes clockwise. (It is clear what length of arcs of $\partial R$ means
and its definition could be easily expanded. For circular arcs it would need more space and 
therefore we replaced circles with boundaries of rectangles.) For any $q\in\C\backslash\{p\}$ we denote by $q_R$ 
the intersection of $\partial R$ 
with the half-line from $p$ through $q$. For any partition $P=(t_0,t_1,\ds,t_k)$, $0=t_0<t_1<\ds<t_k=1$, of $I$
(or of an subinterval $[a,b]\sus I$) we define the {\em winding sum} as
$$
W(P,f,p)=\sum_{j=1}^k\ell(f(t_{j-1})_R,f(t_j)_R)\;.
$$
As before we set $\|P\|=\max_{1\le j\le k}|t_j-t_{j-1}|$. 

\begin{prop}[definition of $w(f,p)$]\tec\label{def_wn}
For every partition $P$ of $I$, we have $W(P,f,p)\in\Z$. If $P_n$ is any sequence of partitions of $I$ with 
$\|P_n\|\to0$ then the integer sequence $W(P_n,f,p)$ is eventually constant. We define, for such sequence $P_n$, 
the winding number of $f$ with respect to $p$ by
$$
w(f,p)=\lim_{n\to\infty}W(P_n,f,p)\;.
$$
\end{prop}
\proof
After cutting $\partial R$ in a point different from all $a_j=f(t_j)_R$ and straightening it we may think of 
$a_j$ as lying in $(0,1)$. Then $\ell(a_{j-1},a_j)=a_j-a_{j-1}+m_j$ for some $m_j\in\{-1,0,1\}$. Thus  
$W(P,f,p)=\sum_{j=1}^k\ell(a_{j-1},a_j)=0+\sum_{j=1}^km_j$ (recall that $a_k=a_0$) is an integer.

It suffices to show that for some $\de>0$, if $P_1,P_2$ are partitions of $I$ with $\|P_1\|,\|P_2\|<\de$ then 
$W(P_1,f,p)=W(P_2,f,p)$. It is easy to see that the telescoping formula $\ell(b_0,b_k)=\sum_{j=1}^k\ell(b_{j-1},b_j)$, 
$b_j\in\partial R$, holds if all $b_j$ lie in an arc of $\partial R$ shorter than $\frac{1}{2}$ (but it fails in 
general). Let first $P_1\sus P_2$, $P_1=(t_0,t_1,\ds,t_k)$, be two partitions of $I$ with $\|P_1\|<\de$. Then there are
partitions $Q_j$ of $[t_{j-1},t_j]$ (using points in $P_2$) such that $W(P_2,f,p)=\sum_{j=1}^kW(Q_j,f,p)$. 
The mapping $I\ni t\mapsto f(t)_R\in\partial R$ is uniformly continuous and therefore if $\de$ is small, each set 
$f([t_{j-1},t_j])_R$ is contained in an arc of $\partial R$ with length less than $\frac{1}{2}$. The telescoping formula
then gives $\ell(f(t_{j-1})_R,f(t_j)_R)=W(Q_j,f,p)$ and so $W(P_1,f,p)=W(P_2,f,p)$. If $P_1,P_2$ are general 
partitions of $I$ with $\|P_1\|,\|P_2\|<\de$, we consider the pairs of partitions $P_1,P_3$ and $P_2,P_3$ where 
$P_3=P_1\cup P_2$ and see again that $W(P_1,f,p)=W(P_3,f,p)=W(P_2,f,p)$. 
\eproof

We prove that the $w(f,p)$ we have just defined has the properties stated in Proposition~\ref{wind_numb} 
(at the end of Comments we sketch another way how to define $w(f,p)$). Integrality has been already proven. For any  
constant loop $f$, every winding sum is a sum of zeros, so it is zero and $w(f,p)=0$. Let $p$ be the center of a 
perimeter $1$ rectangle $R$ and $f:\;I\to\partial R$ be a positive circuit. Then for any partition 
$P=(t_0,t_1,\ds, t_k)$ of $I$ with small enough $\|P\|$, for the points $a_j=f(t_j)_R=f(t_j)$, $j=0,1,\ds,k$, each 
arc of $\partial R$ obtained by going from $a_{j-1}$ counter-clockwisely to $a_j$ has length less than $\frac{1}{2}$.
Thus $W(P,f,p)=\sum_{j=1}^k\ell(a_{j-1},a_j)$ is just the sum of these lengths, which is $p(R)=1$, and $w(f,p)=1$. 
General rectangles and general interior points reduce to this case with the help of properties 2 and 3, which we prove 
in a moment. Thus property 1 is proven, modulo properties 2 and 3.

Properties 2, 3, and 4 follow from the fact that for any $\de>0$ the mapping $F:\;\{(p,q)\in\C^2\;|\;|p-q|>\de\}\to\C$, 
$F(p,q)=q_{R_p}=\partial R_p\cap\{t(q-p)\;|\;t\ge0\}$, is uniformly continuous; note that in this case the domain 
of $F$ is not compact. Suppose that $p_1,p_2\in K$ are 
two points of the same component $K$ of $\C\backslash f(I)$ where $f$ is a loop. We take a path $g:\;I\to K$
joining them in $K$, $g(0)=p_1$ and $g(1)=p_2$. By the uniform continuity of $F$ the mapping $I\times I\ni (u,t)
\mapsto f(t)_{R_{g(u)}}$ is uniformly continuous because $g(I)$ and $f(I)$ have positive distance (as disjoint 
compact sets). So there is a $\de>0$ such that for any fixed partition $P$ of $I$ with $\|P\|<\de$ the mapping 
$I\ni u\mapsto W(P,f,g(u))\in\Z$ is continuous and hence constant (small $\de$ ensures that for every $u\in I$ 
all arcs $a_{j-1}$-$a_j$ of $\partial R_{g(u)}$ involved in $W(P,f,g(u))$ are shorter than, say, 
$\frac{1}{4}<\frac{1}{2}$ and $\ell(\cdot,\cdot)$ 
varies continuously). So $W(P,f,p_1)=W(P,f,p_2)$ and hence $w(f,p_1)=w(f,p_2)$, which is property 2. 
Suppose that $h:\;I\times I\to\C$ is continuous with $h(0,u)=h(1,u)$ for every $u\in I$ and $p\in\C\backslash h(I\times I)$ 
is a point. As before, 
by the uniform continuity of $F$ (and compactness of $h(I\times I)$) there is a $\de>0$ such that for any 
fixed partition $P$ of $I$ with $\|P\|<\de$ the mapping $I\ni u\mapsto W(P,h(t,u),p)\in\Z$  is continuous and 
hence constant. So $W(P,h(t,0),p)=W(P,h(t,1),p)$ and hence $w(h(t,0),p)=w(h(t,1),p)$, which is property 3. 
Property 4 is of course a special case of property 3, for the homotopy $h(t,u)=f(us+t)$.

Let $h=f*g$ be a product of loops and $p\in\C\backslash h(I)$ be a point. Property 5 follows at once from 
the identity $W(PQ,h,p)=W(P,f,p)+W(Q,g,p)$ where $P=(t_0,\ds,t_k)$ and $Q=(u_0,\ds,u_l)$ are arbitrary partitions 
of $I$ and $PQ=(v_0,\ds,v_{k+l})$ is the concatenated partition of $I$ given by $v_i=t_i/2$ for $0\le i\le k$ 
and $v_i=(u_{i-k}+1)/2$ for $k\le i\le k+l$.

\bigskip\noindent
Only now is the second proof of Theorem~\ref{zakl_veta} complete.

\subsection{Comments}

Continuing in this way, Connell \cite{conn65} proves the existence of $f''$ and thus of all derivatives $f^{(n)}$, 
and Porcelli and Connell \cite{porc_conn} deduce the full Taylor expansion of $f$. 

The previous integration-free proof of Theorem~\ref{zakl_veta} roughly follows the proof of Connell 
\cite{conn65} (for example, part 1 of Proposition~\ref{cor_w_num} is \cite[Lemma 1]{conn65}, part 3 of 
Proposition~\ref{zero_deri} is \cite[Lemma 4]{conn65}, and  
Proposition~\ref{end_2nd_proof} is \cite[Theorem 2]{conn65}) who acknowledges that he borrows from Whyburn \cite{whyb}, 
but we filled in proofs and precise statements of many steps and results that
Connell omits or takes for granted, such as the definition and properties of $w(f,p)$, the fact that an open 
set does not have measure $0$, and others. Also, familiar auxiliary results like 
that zero derivative implies constant function have to be established without integration. 
Now, when we filled in all details and omitted steps, we can fairly compare lengths of both proofs, the integration 
one versus the topological one. 
Topology indeed cannot compete with integration as the score settles at less than $4$ pages to 
$8\frac{1}{2}$ pages. But 
several techniques employed in the topological proof, such as the winding number or K\"onig's lemma, are quite nice
and open connections to other areas of mathematics.

Plunkett \cite{plun} gave an integration-free proof of continuity of complex derivative in 1959 but was unaware 
that this had been already achieved by Adel'son-Vel'skii and Kronrod \cite{adel_kron} in the USSR 15 years earlier, 
as pointed out by Lorent \cite{lore}; we will discuss their proof in part 2. Topological techniques developed
by Whyburn \cite{whyb} were used by him \cite{whyb61}, Connell \cite{conn61,conn65}, Connell and Porcelli 
\cite{conn_porc,porc_conn}, Plunkett \cite{plun,plun60} and Read \cite{read} to build the theory of analytic functions in 
an integration-free way on a topological base. Eggleston and Ursell \cite{eggl_urse} also deduced properties 
of analytic functions like the maximum modulus principle by the winding number. 

K\"onig's lemma, which we used in the proof of part 2 of Proposition~\ref{konig}, is due to K\H onig \cite{koni} 
in 1927 and asserts that every infinite but finitely branching rooted tree has an infinite branch. It appears implicitly 
in the proof of Cauchy's theorem in Titchmarsh \cite{titc} and elsewhere, indeed already in Goursat 
\cite[the proof of Lemme on p. 15]{gour} in 1900. 

For an exposition, connections, and applications of the winding 
number see the book of Roe \cite{roe}. A spectacular application of the winding number, more precisely of 
its discrete version, is in the rigorous proof by Duminil-Copin and Smirnov \cite{dumi_smir} of the conjecture 
due to Nienhuis \cite{nien} that the numbers $a_n$ of $n$-step self-avoiding walks from a fixed point 
in the honeycomb lattice (the graph of the tiling 
of the plane by regular hexagons) have growth 
$$
\lim_{n\to\infty}a_n^{1/n}=\sqrt{2+\sqrt{2}}\;.
$$ 
For further applications of the winding number method to counting self-avoiding walks see Beaton, Bousquet-M\'elou, 
de Gier, Duminil-Copin and Guttmann \cite{beat_al} and Glazman \cite{glaz}. See Kahane \cite{kaha} and Bourgain 
and Kozma \cite{bour_kozm} for expressing $w(f,p)$ in terms of Fourier coefficients. One often defines the winding 
number by the integral 
$$
w(f,p)=\frac{1}{2\pi i}\int_f\frac{dz}{z-p} 
$$
but this is not available when integration is not allowed. 
We admit that our definition in Proposition~\ref{def_wn} is a thinly disguised integration.
One can avoid integration by another standard definition of $w(f,p)$ which we indicate now. 

For two sets $A,B\sus\R$ we define their distance by 
$$
d(A,B)=\inf_{a\in A,b\in B}|a-b|\;.
$$
For a fixed point $p\in\C$ we take the perimeter $1$ square $R=R_p$ centered at $p$ and for any $z\in\C\backslash\{p\}$
denote by $l(z)\in[0,1)$ the length of the counter-clockwise $a$-$z_R$ arc of $\partial R$ (recall that $a$ is the 
lower left corner of $R$ and $z_R$ is the intersection of $\partial R$ with the ray from $p$ through $z$). 
The coordinate mapping $\C\backslash\{p\}\ni z\mapsto l(z)\in[0,1)$ is discontinuous, which we fix
by replacing it with the mapping $z\mapsto P(z):=\Z+l(z)=\{n+l(z)\;|\;n\in\Z\}$. This mapping is continuous to 
$d(\cdot,\cdot)$; $d(\Z+a,\Z+b)$, $a,b\in\R$, in fact equals the distance of $a-b$ to the nearest integer. 
One can show that 
for any path $f:\;I\to\C\backslash\{p\}$ there exists a continuous function $g:\;I\to\R$ such that $g(t)\in P(f(t))$ 
for every $t\in I$, and that $g$ is unique up to an integral shift. Thus $g(1)-g(0)$ is independent of $g$, and if $f$ 
is a loop then $g(1)-g(0)\in\Z$. One can prove Proposition~\ref{wind_numb} by defining $w(f,p)$ as the difference 
$g(1)-g(0)$. This construction of $w(f,p)$ is usually stated instead of $P(z)$ in terms of the argument of the 
complex logarithm $\log z$, see for example Whyburn \cite{whyb61}.

\section{Conclusion and outlook}

We presented two proofs of Theorem~\ref{zakl_veta}. The third proof by Adel'son-Vel'skii and Kronrod 
\cite{adel_kron1,adel_kron2,adel_kron} will follow, let us hope, in part 2. The amount of material to be filled 
in to make it rigorous and self-contained is considerably larger compared to the second proof.

\bigskip\noindent
Martin Klazar\\
Department of Applied Mathematics (KAM)\\
Faculty of Mathematics and Physics\\
Charles University\\ 
Malostransk\'e n\'am. 25\\ 
118 00 Praha 1\\
Czechia\\
{\tt klazar@kam.mff.cuni.cz}


\begin{thebibliography}{XX}
\bibitem{maov} Is there an integration free proof (or heuristic) that once differentiable implies twice
differentiable for complex functions?\\
{\tt mathoverflow.net/questions/63205}
\end{thebibliography}
\end{document}